\documentclass[11pt]{article}
\usepackage[leqno]{amsmath} 
\usepackage{amssymb,amscd,amsthm}
\input xy \xyoption{all}
\CompileMatrices

\def\Ch{{\cal C}}
\def\Ah{{\cal A}}
\def\Oh{{\cal O}}

\def\Eh{{\cal E}}
\def\Dh{{\mathfrak{D}}}
\def\Hh{{\cal H}}
\def\Ih{{\cal I}}
\def\Lh{{\cal L}}
\def\Mh{{\cal M}}
\def\Nh{{\cal N}}

\def\Th{{\cal T}}

\def\Uhat{{\hat{\cal U}}}
\def\Lh{{\cal L}}

\def\Z{\mathbb{Z}}
\def\Q{\mathbb{Q}}
\def\Qbar{\bar{\mathbb{Q}}}
\def\C{\mathbb{C}}
\def\R{\mathbb{R}}
\def\P{\mathbb{P}}
\def\ffr{\mathfrak{f}}
\def\bfr{\mathfrak{b}}
\def\afr{\mathfrak{a}}

\def\hfr{\mathfrak{H}^g_\pm}

\def\dS{{\partial S}}
\def\Sbar{{\overline{S}}}

\def\tildeh{{\widetilde{h}}}

\def\isom{\cong}
\def\verk{\circ}
\def\prolim{\varprojlim}

\def\Hom{\operatorname{Hom}}
\def\Eis{\operatorname{Eis}}

\def\Ext{\operatorname{Ext}}

\def\Gl{\operatorname{Gl}}

\def\Gm{{\mathbb{G}_m}}

\def\Tr{\operatorname{Tr}}

\def\id{\operatorname{id}}
\def\pol{\operatorname{pol}}
\def\Log{\operatorname{Log}}

\def\Sym{\operatorname{Sym}}

\def\vol{\operatorname{vol}}
\def\rat{\operatorname{rat}}
\def\res{\operatorname{res}}
\def\Res{\operatorname{Res}}

\def\tors{\operatorname{tors}}

\def\Ind{\operatorname{Ind}}

\def\Spec{\operatorname{Spec}}

\newcommand{\bew}{\begin{proof}}
\newcommand{\bewende}{\end{proof}}
\newtheorem{lemma}{Lemma}[subsection]
\newtheorem{prop}[lemma]{Proposition}
\newtheorem{thm}[lemma]{Theorem}
\newtheorem{defn}[lemma]{Definition}
\newtheorem{cor}[lemma]{Corollary}

\newcommand{\rem}{\noindent{\bf Remark:\ }}

\title{\vskip -25mm Degeneration of polylogarithms and special values
of L-functions for totally real fields}   
\author{Guido Kings}
\date{December 7, 2005}

\begin{document}

\maketitle


\section*{Introduction} 
The polylogarithm is a very powerful tool in studying special values
of $L$-functions and subject to many conjectures. Most notably,
the Zagier conjecture claims that all values of $L$-functions
of number fields can be described by polylogarithms. The interpretation
of the polylogarithm functions in terms of periods of variations of
Hodge structures has lead to a motivic theory of the polylog and to
generalizations as the elliptic polylog by Beilinson and Levin.
Building on this work, Wildeshaus has defined polylogarithms in
a more general context and in particular for abelian schemes. 

Not very much is known about the extension classes arising from these
``abelian polylogarithms''. In an earlier paper \cite{K} we were able to show that the abelian polylogarithm, as defined by Wildeshaus, is indeed of motivic origin, i.e., is in the
image of the regulator from $K$-theory. 

It was Levin in \cite{L}, who started to investigate certain "polylogarithmic currents" 
on abelian schemes, which are related to the construction by Wildeshaus. 
Very recently, Blotti\`ere could show in his thesis \cite{B}
that these currents actually represent the polylogarithmic extension in the category 
of Hodge modules. Furthermore, specializing to the case of 
Hilbert modular varieties, he computed the residue of the associated
Eisenstein classes, which are just the pull-back of the polylog along torsion sections
of the universal abelian variety. This residue is given in terms
of (critical) special values of $L$-functions of the totally real field,
which defines the Hilbert modular variety. His computation uses the explicit
description of the polylog in terms of the currents constructed by Levin \cite{L}.
%

In this paper we will, following and extending ideas from the case 
of the elliptic polylog treated in \cite{HK} 
(which is in turn inspired by \cite{BL}), present a completely different 
approach to this residue computation, which avoids computations as much as possible and
relates the polylog on the Hilbert modular variety to the classes constructed 
by Nori and Sczech. Instead of computing directly
the degeneration on the base (as in the approach by Blotti\`ere), 
we work with the polylogarithm, which lives on
the universal abelian scheme and use its good functorial properties to
compute the degeneration. Contrary to earlier approaches to
the degeneration, which work in the algebraic category of schemes, we
view the problem as of topological
nature and work entirely in the topological category. For this
the notion of "topological polylogarithm", which is defined for an
arbitrary real torus, is essential. We think that this unusual approach is of
independent interest.

The idea to study the polylogarithm on the moduli scheme of elliptic curves 
via its degeneration at the boudary is already prominent in \cite{BL}, where it
is shown that the elliptic polylog degenerates into the cyclotomic polylog. 
These ideas were developed further by Wildeshaus \cite{W2} in the context of toroidal
compactification of moduli schemes
of abelian varieties. He could show, that also in this general context
the polylogarithm is stable under degeneration. In \cite{HK} this degeneration
principle was exploited for the moduli space of elliptic curves and used to relate the elliptic polylogarithm to the critical and non-critical values of Dirichlet L-functions.

To describe the theorem more precisely, consider the specialization
of the polylog, which gives Eisenstein classes (say in the category 
of mixed \'etale sheaves to fix ideas) 
\[
\Eis^k(\alpha)\in \Ext^{2g-1}_S(\Q_\ell, \Sym^k\Hh(g)),
\]
where $S$ is the Hilbert modular variety of dimension $g$ and $\Hh$ is 
the locally constant sheaf of relative Tate-modules
of the universal abelian scheme. 
Let $j:S\to\Sbar$ be the Baily-Borel compactification of $S$
and  $i:\dS:= \Sbar\setminus S\to \Sbar$ the
inclusion of the cusps. The degeneration or residue 
map (see \ref{deg2} for the precise definition) is then
\[
\res:\Ext^{2g-1}_S(\Q_\ell, \Sym^k\Hh(g))\to \Hom_\dS(\Q_\ell,\Q_\ell).
\]
The target of this map is sitting inside a sum of copies of $\Q_\ell$ and
the main result of this paper \ref{maintheorem} describes 
$\res(\Eis^k(\alpha))$ in
terms of special values of (partial) $L$-functions of the totally real
field defining $S$.

There is a very interesting question raised by the results in this paper. 
In \cite{HK} we were able to construct extension
classes related to non-critical values of Dirichlet-$L$-functions, if the
residue map was zero on the specialization of the polylog. Is there an
analogous result here? 

The paper is organized as follows: In the first section we review the
definition of the Hilbert modular variety, define the residue or
degeneration map and formulate our main theorem. 
The second section reviews the theory of the polylog and the Eisenstein
classes
emphasing the topological situation, which is not extensively covered 
in the literature. In the third section we give the proof of the
main theorem. 

It is a great pleasure to thank David Blotti\`ere for a series of interesting 
discussions during his stay in Regensburg, which stimulated the research in this paper. 
Moreover, I like to thank Sascha Beilinson for making available some time ago his notes
about his and A. Levin's interpretation of Nori's work. Finally, I thank J\"org 
Wildeshaus for helpful comments.

\section{Polylogarithms and degeneration}

We review the definition of  a Hilbert modular variety to fix 
notations and
pose the problem of computing the degeneration of the
specializations of the polylogarithm at the boundary. The main theorem
describes this residue in terms of special values of $L$-functions.

\subsection{Notation}\label{notation}
As in \cite{BL} we deal with three different types of sheaves simultaneously.
Let $X/k$ be a variety and $L$ a coefficient ring for our sheaf theory, 
then we consider
\begin{itemize}
\item[i)] $k=\C$ the usual topology on $X(\C)$ and $L$ any commutative ring
\item[ii)] $k=\R$ or $\C$ and $L=\Q$ or $\R$ and we work with the category of
mixed Hodge modules 
\item[iii)] $k=\Q$ and $L=\Z/l^r\Z,\Z_l$ or $\Q_l$ 
and we work with the category of {\'e}tale sheaves
\end{itemize}

\subsection{Hilbert modular varieties}
We recall the definition of Hilbert modular varieties following 
Rapoport \cite{R}. To avoid all technicalities, we will only
consider the moduli scheme over $\Q$. The theory works over more general
bases schemes without any modification.

Let $F$ be a totally real field, $g:=[F:\Q]$, $\Oh$ the ring of integers,
$\Dh^{-1}$ the inverse different and $d_F$ its discriminant. 
Fix an integer  $n\geq 3$.
We consider the functor, which associates to a scheme $T$ over 
$\Spec\Q$
the isomorphism classes of triples $(A,\alpha,\lambda)$, where $A/T$ is
an abelian scheme of dimension $g$, with real multiplication by 
$\Oh$, $\alpha:\underline{\Hom}_{et,\Oh,sym}(\Ah, \Ah^*)\to \Dh^{-1} $ 
is a $\Dh^{-1}$-polarization in the sense of \cite{R} 1.19, i.e., an 
$\Oh$-module isomorphism respecting
the positivity of the totally positive elements in $\Dh^{-1}\subset F$,  
and $\lambda:A[n]\isom (\Oh/n\Oh)^2$ is a level $n$ structure satisfying 
the compatibility of \cite{R} 1.21. For $n\geq 3$ this functor 
is represented by a
smooth scheme $S:=S^{\Dh^{-1}}_n$ of finite type over $\Spec \Q$ .
Let 
\[
\Ah\xrightarrow{\pi}S
\]
be the universal abelian scheme over $S$. In any of the three categories
of sheaves i)-iii) from \ref{notation} we let
\[
\Hh:= \underline{\Hom}_S(R^1\pi_*L,L)
\]
the first homology of $\Ah/S$. In the {\'e}tale case and $L=\Z_\ell$, 
the fiber of $\Hh$ at a point is the Tate module of the abelian variety over
that point.

\subsection{Transcendental description}
For the later computation we need a description in group theoretical
terms of the complex points $S(\C)$ and of $\Hh$.

Define a group scheme $G/\Spec \Z$ by the Cartesian diagram
\[
\begin{CD} G@>>>\Res_{\Oh/\Z}\Gl_2\\
@VVV@VV\det V\\
\Gm@>>> \Res_{\Oh/\Z}\Gm
\end{CD}
\]
and let 
\[
\hfr:= \{\tau\in F\otimes\C|Im\; \tau \mbox{ totally positive or totally negative}\}.
\] 
Then $\left(\begin{array}{cc}a&b \\ c&d\end{array}\right)\in G(\R)$ acts
on $\hfr$ by the usual formula
\[
\left(\begin{array}{cc}a&b \\ c&d\end{array}\right)\tau=
\frac{a\tau+b}{c\tau+d}
\]
and the stabilizer of $1\otimes i\in \hfr$ is 
\[
K_\infty:= (F\otimes \C)^*\cap G(\R),
\]
so that 
\[
\hfr\isom G(\R)/K_\infty.
\]
With this notation one has
\[
S(\C)=G(\Z)\backslash (\hfr\times G(\Z/n\Z)).
\]
On $S(\C)$ acts $G(\Z/n\Z)$ by right multiplication.
The determinant $\det:G\to \Gm$ induces 
\[
S(\C)\to \Gm(\Z/n\Z)
\]
and the fibers are the connected components. Define a subgroup 
$D\subset G$ isomorphic to $\Gm$ by 
$D:=\{\bigl(\begin{smallmatrix} a&0\\ 0&1
\end{smallmatrix}\bigr)\in G:a\in \Gm \bigr)$. This gives a section of
$\det$.
Then the action of $D(\Z/n\Z)$ by right multiplication  is transitive
on the set of connected components. 

The embedding $G(\Z)\subset \Gl_2(\Oh)$ defines an action of
$G(\Z)$ on $\Oh^{\oplus 2}$ and in the topological 
realization the local system $\Hh$ is given by the quotient
\[
G(\Z)\backslash (\hfr\times\Oh^{\oplus 2} \times  G(\Z/n\Z)).
\]
In particular, as a family of real ($2g$-dimensional) tori, the complex
points $\Ah(\C)$ of the universal
abelian scheme can be written as 
\[
G(\Z)\backslash (\hfr\times (F\otimes\R/\Oh)^{\oplus 2}\times  G(\Z/n\Z))
\]
and the level $n$ structure is given by the subgroup 
\[
(\frac{1}{n}\Oh/\Oh)^{\oplus 2}\subset (F\otimes\R/\Oh)^{\oplus 2}.
\]
The $\Oh$-multiplication on $\Ah(\C)$ is in this description given
by the natural $\Oh$-module structure on $F\otimes\R$.
\subsection{Transcendental description of the cusps}
The following description of the boundary cohomology is inspired by
\cite{Ha}.
Let $B\subset G$ the subgroup of upper triangular matrices,
$T\subset B$ its maximal torus and $N\subset B$ its unipotent radical.
We have an exact sequence
\[
1\to N\to B\xrightarrow{q} T\to 1.
\]
We denote by $G^1$, $B^1$ and $T^1$ the subgroups of  determinant $1$.
Let $K^B_\infty:=B(\R)\cap K_\infty$, then the
Cartan decomposition shows that $\hfr=B(\R)/K^B_\infty$.
A pointed neighborhood of the set of all cusps is given by
\begin{equation}
\widetilde{S}_B:=B(\Z)\backslash\bigl( B(\R)/K^B_\infty \times G(\Z/n\Z)\bigr).
\end{equation}
In particular, the set of cusps is
\begin{equation}\label{cusps}
\dS(\C)=B^1(\Z)\backslash G(\Z/n\Z).
\end{equation}
The fibres of the map $\dS(\C)\to \Gm(\Z/n\Z)$ induced by the determinant are 
\begin{equation}\label{cuspsasp1}
B^1(\Z)\backslash G^1(\Z/n\Z)\isom \Gamma_G\backslash \P^1(\Oh),
\end{equation}
where $\Gamma_G:=\ker(G^1(\Z)\to G^1(\Z/n\Z))$. In particular, we can
think of a cusp represented by  $h\in  G^1(\Z/n\Z)$ as a rank $1$ 
$\Oh$-module $\bfr_h$, which is a quotient
\begin{equation}\label{phdefn}
\Oh^2\xrightarrow{p_h}\bfr_h,
\end{equation}
together with a level structure, i.e., a basis $h\in G^1(\Z/n\Z)$. Explicitly,
the fractional ideal $\bfr_h$ is generated by any representatives $u,v\in\Oh$
of the second row of $h$.

On $\widetilde{S}_B$ acts $G(\Z/n\Z)$ by multiplication from the right.
This action is transitive on the connected 
components of $\widetilde{S}_B$. 
Define
\begin{equation}\label{cuspneighbor}
{S}_B:=B(\Z)\backslash\bigl( B(\R)/K^B_\infty \times B(\Z/n\Z)\bigr),
\end{equation}
then $S_B\subset \widetilde{S}_B$ is a union of connected components of
$\widetilde{S}_B$.
Let $K^T_\infty$ (respectively $T(\Z)$) be the image of
$K^B_\infty$ (respectively $B(\Z)$) under $q:B(\R)\to T(\R)$. Define
\begin{equation}
S_T:= T(\Z)\backslash \bigl(T(\R)/K^T_\infty\times T(\Z/n\Z)\bigr),
\end{equation}
then the map $q:B\to T$ induces a fibration
\begin{equation}
q:S_B\to S_T,
\end{equation}
whose fibers are $N(\Z)\backslash \bigl(N(\R)\times N(\Z/n\Z)\bigr)$ with 
$N(\Z):=B(\Z)\cap N(\R)$. Denote by 
\begin{equation}\label{udefn}
u:S_T\to pt
\end{equation}
the structure map to a point.
For the study of the degeneration, one considers
the diagram
\begin{equation}\label{comdiagr}
\begin{CD} S_B@>q>> S_T@>u >>pt\\
@VVV\\
S
\end{CD}
\end{equation}
In fact we are interested in the cohomology of certain local
systems on these topological spaces. For the computations it
is convenient to replace $S_B$ and $S_T$ by homotopy equivalent spaces
as follows. 

Define $K^{T^1}_\infty:=K^T_\infty\cap T^1(\R)$. Then the inclusion
induces an isomorphism
\[
S^1_T:=T^1(\Z)\backslash \bigl( T^1(\R)/K^{T^1}_\infty\times T(\Z/n\Z)\bigr)\isom S_T.
\]
The map 
$a\mapsto \bigl(\begin{smallmatrix} a&0\\ 0& a^{-1}\end{smallmatrix}\bigr)$
defines isomorphisms $ (F\otimes \R)^*\isom T^1(\R)$ and
$\Oh^*\isom T^1(\Z)$. Note that $K^{T^1}_\infty\subset (F\otimes \R)^*$
is identified with the two torsion subgroup in $(F\otimes \R)^*$ and that
$K^{T^1}_\infty\isom (\Z/2\Z)^g$ 
permutes the set of connected components of $T^1(\R)$.

\begin{lemma}\label{s1defn} Let $(F\otimes \R)^1$ be the subgroup of $(F\otimes \R)^*$
of elements of norm $1$ and $\Oh^{*,1}=\Oh^*\cap (F\otimes \R)^1$. Then 
\[
S^1_T:=\Oh^{*,1}\backslash 
\bigl((F\otimes \R)^1/K^{T^1}_\infty\cap(F\otimes \R)^1 \times T(\Z/n\Z)\bigr)
\]
is homotopy equivalent to $S_T$. Moreover, the inclusion of
the totally positive elements $(F\otimes \R)^1_+$
into $(F\otimes \R)^1$ provides an identification
\[
(F\otimes \R)^1_+\isom(F\otimes \R)^1/K^{T^1}_\infty\cap(F\otimes \R)^1 .
\]
\end{lemma}
\bew 
The exact sequence
\[
0\to (F\otimes \R)^1\to (F\otimes \R)^*\to \R^*\to 0
\]
together with the fact that $K^{T^1}_\infty$ is the two torsion in
$(F\otimes \R)^*$ allows to identify
\[
T^1(\R)/K^{T^1}_\infty\isom \bigl((F\otimes \R)^1/K^{T^1}_\infty\cap(F\otimes \R)^1\bigr)\times \R_{>0}.
\]
The last identity is clear.
\bewende
We define $S^1_B$ to be the inverse image of $S^1_T$ under $q$, so that
we have a Cartesian diagram
\begin{equation}\label{eq:basedegener}
\begin{CD}
S^1_B@>q>> S^1_T\\
@VVV@VVV\\
S_B@>q>>S_T.
\end{CD}
\end{equation}
Over $S_B$ the representation $\Oh^2$ has a filtration
\begin{align}
  \label{eq:localsystemdegener}
0\to \Oh\to \Oh^2\xrightarrow{p} \Oh\to 0,
\end{align}
where the first map
sends $a\in \Oh$ to the vector ${a\choose 0}$ and the second map
is ${a\choose b}\mapsto b$.
This induces a filtration on the
local system $\Hh$ 
\begin{equation}\label{locexseq}
0\to \Nh\to \Hh\to \Mh\to 0,
\end{equation}
where $\Nh$ and $\Mh$ are the associated local systems.
In particular, over $S^1_B$ one has a filtration
of topological tori
\begin{equation}
  \label{eq:toridegener}
  0\to\Th_\Nh\to \Ah(\C)\xrightarrow{p} \Th_\Mh\to 0,
\end{equation}
where $\Th_\Nh:= \Nh\otimes\R/\Z$ and $\Th_\Mh:=\Mh\otimes\R/\Z$.
By definition of $N$ the fibration in 
(\ref{eq:toridegener}) and (\ref{eq:basedegener}) are compatible, i.e., one
has a commutative diagram
\begin{equation}\label{degenerationdiagr}
\begin{CD}
\Ah(\C)@>p>>\Th_\Mh\\
@V\pi VV@VV\pi_\Mh V\\
S^1_B @>q>> S^1_T.
\end{CD}
\end{equation}

\subsection{The degeneration map}
In this section  we explain the degeneration problem we want to consider.

The polylogarithm on $\pi:\Ah\to S$ defines for certain linear
combinations $\alpha$ of torsion sections  of $\Ah$ an extension class
\begin{equation}\label{pol}
\Eis^k(\alpha)\in \Ext^{2g-1}_{?,S}(L, \Sym^k\Hh(g)),
\end{equation}
where $?$ can be $MHM,et,top$.
The construction of this class will be given in section \ref{polsection}
definition \ref{eisdefn}. 

Let $\Sbar$ be the Baily-Borel compactification of $S$. Denote
by $\dS:=\Sbar\setminus S$ the set of cusps. We get
\[
\dS\xrightarrow{i}\Sbar\xleftarrow{j}S.
\]
The adjunction map together with the edge morphism in the Leray spectral
sequence for $Rj_*$ gives
\begin{equation}\label{deg1}
\begin{CD}
\Ext^{2g-1}_S(L, \Sym^k\Hh(g))@>>> \Ext^{2g-1}_{\dS}(L, i^*Rj_*\Sym^k\Hh(g))\\
&\searrow&@VVV\\
&& \Hom_{\dS}(L, i^*R^{2g-1}j_*\Sym^k\Hh(g)).
\end{CD}
\end{equation}
There are several possibilities to compute $i^*R^{2g-1}j_*\Sym^k\Hh(g)$.

\begin{thm}\label{degenercompu} Assume that $\Q\subset L$. Then, in any of
the categories $MHM,et,top$,
there is a canonical isomorphism
\[
i^*R^{2g-1}j_*\Sym^k\Hh(g)\isom L,
\]
where $L$ has the trivial Hodge structure (resp. the trivial Galois action).
\end{thm}
\rem J\"org Wildeshaus has pointed out that the determination of the weight on
the right hand side is not necessary for  our main result, but follows
from it. In fact, our main result gives non-zero classes in
\[
\Hom_{\dS}(L, i^*R^{2g-1}j_*\Sym^k\Hh(g)),
\]
so that the rank one sheaf
$i^*R^{2g-1}j_*\Sym^k\Hh(g)$ has to be of weight zero.\\

Using this identification we define the residue or degeneration map:
\begin{defn}\label{deg2} 
The map from (\ref{deg1}) together with the identification
of \ref{degenercompu} define
the {\em residue map}
\[
\res:\Ext^{2g-1}_S(L, \Sym^k\Hh(g))\to\Hom_{\dS}(L,L).
\]
The residue map is equivariant for the $G(\Z/n\Z)$ action on
both sides.
\end{defn}

\bew (of theorem \ref{degenercompu}). 
In the case of Hodge modules we use theorem 2.9. in Burgos-Wildeshaus \cite{BW}
and in the \'etale case we use theorem 5.3.1 in \cite{Pi}.
Roughly speaking, both results   
asserts that the higher direct image can be calculated using
group cohomology and the ``canonical construction'', which associates 
to a representation of the group defining the Shimura variety a Hodge 
module resp. an \'etale sheaf. 

More precisely, from a topological point of view, the monodromy at the cusps
is exactly the cohomology of $\widetilde{S}_B$. One has
\[
H^{2g-1}(\widetilde{S}_B, \Sym^k\Hh(g))\isom\Ind_{B(\Z/n\Z)}^{G(\Z/n\Z)}
H^{2g-1}(S_B,\Sym^k\Hh(g))
\]
and 
\[
H^{2g-1}(S_B,\Sym^k\Hh(g))\isom \bigoplus_{r+s=2g-1}
H^r(S_T,R^sq_*\Sym^k\Hh(g))).
\]
As the cohomological dimension of $\Gamma_T$ is $g-1$ and that of $\Gamma_N$ is
$g$, one has in fact
\[
H^{2g-1}(S_B,\Sym^k\Hh(g))\isom H^{g-1}(S_T,R^gq_*\Sym^k\Hh(g))).
\]
The exact sequence 
\[
0\to \Oh\to \Oh^2\xrightarrow{p} \Oh\to 0
\]
from (\ref{eq:localsystemdegener}) shows that $R^gq_*\Sym^k\Hh(g)$ 
can be identified via $p$ with $\Sym^k\Oh\otimes L$ 
with the induced $T(\Z)$ action,
which maps $\bigl(\begin{smallmatrix} a&0\\ 0& d\end{smallmatrix}\bigr)$
to $d^{k}$. To compute the coinvariants, extend the coefficients to $\R$, 
so that 
\[
\Oh\otimes \R\isom \bigoplus_{\tau:F\to\R} \R
\]
and $\bigl(\begin{smallmatrix} a&0\\ 0& d\end{smallmatrix}\bigr)\in T(\Z)$
acts via $\tau(d)$ on the  component indexed by $\tau$. 
Thus $\Sym^k\Oh\otimes L$ 
can only have a trivial quotient, if $k\equiv 0\mbox{ mod }g$ and on
this one dimensional quotient the action is by the norm map $T(\Z)\to \pm 1$.
One gets:
\[
H^{g-1}(S_T, \Sym^k\Oh\otimes L)\isom \left\{\begin{array}{cc} L&\mbox{ if } k
\equiv 0
\mbox{ mod }g\\
0&\mbox{ else}
\end{array}\right.
\]
The above mentioned theorems imply that this topological computation
gives also the result in the categories $MHM, et, top$. 
The Hodge structure on $H^{g-1}(S_T, \Sym^k\Oh\otimes L) $ is the trivial one,
as one sees from the explicit description of the 
action of $T$ and the fact that the action of the
Deligne torus $\mathbb{S}$, which defines the weight, is induced from
the embedding $x\mapsto \bigl(\begin{smallmatrix} x&0\\ 0& 1\end{smallmatrix}\bigr)$.,
hence is trivial. The same remark and proposition
5.5.4. in \cite{Pi} show that the weight is also zero in the \'etale case.
\bewende

\subsection{Partial zeta functions of totally real fields}

Let $\bfr, \ffr$ be relatively prime integral ideals of $\Oh$, 
$\epsilon: (\R\otimes F)^*\to \{\pm 1\}$ a sign character.
This is a product of characters $\epsilon_\tau:\R^*\to \{\pm 1\}$
for all embeddings $\tau:F\to \R$. Denote by
$|\epsilon |$ the number of non-trivial $\epsilon_\tau$ which occur
in this product decomposition of $\epsilon$. Moreover let 
$x\in \Oh$ such that $x\not\equiv 0 \mbox{ mod }\bfr^{-1}\ffr$ and
$\Oh_{\ffr}^*:=\{ a\in \Oh|a \mbox{ totally positive and }a\equiv 1 \mbox{ mod }\ffr\}$.
Define
\begin{equation}
F(\bfr, \ffr,\epsilon,x, s):= \sum_{\nu\in (x+\ffr\bfr^{-1})/\Oh_{\ffr}^*}
\frac{\epsilon(\nu)}{|N(\nu)|^s}
\end{equation}
for $Re\; s>1$. Here $N$ is the norm. 
On the other hand let $\Tr:F\to \Q$ be the trace map
and define
\begin{equation}
L(\bfr, \ffr,\epsilon,x, s):= \sum_{\lambda\in \bfr(\ffr\Dh)^{-1}/\Oh_{\ffr}^*}
\frac{\epsilon(\lambda)e^{2\pi i\Tr (x\lambda)}}{|N(\nu)|^s}
\end{equation}
These two $L$-functions are related by a functional equation. To formulate
it we introduce the $\Gamma$-factor
\[
\Gamma_\epsilon(s):=\pi^{-\frac{1}{2}(sg+|\epsilon|)}
\Gamma\left(\frac{s+1}{2}\right)^{|\epsilon|}
\Gamma\left(\frac{s}{2}\right)^{g-|\epsilon|}.
\]
The functional equation follows directly with Hecke's method for
Gr\"ossencharacters and was first mentioned by Siegel: 
\begin{prop}[cf.\cite{Si} Formel (10)] The functional equation reads:
\[
\Gamma_\epsilon(1-s)F(\bfr, \ffr,\epsilon,x,1- s)=
i^{-|\epsilon|}|d_F|^{-\frac{1}{2}}N(\ffr^{-1}\bfr)\Gamma_\epsilon(s)
L(\bfr, \ffr,\epsilon,x, s),
\]
where $d_F$ is the discriminant of $F/ \Q$.
\end{prop}
The functional equation shows that $F(\bfr, \ffr,\epsilon,x,1- k)$
can be non-zero for $k=1,2,\ldots$ 
only if $|\epsilon|$ is either $g$ or $0$. Let us introduce
\[
\zeta(\bfr,\ffr,x,s):= \sum_{\nu\in (x+\ffr\bfr^{-1})/\Oh_{\ffr}^*}
\frac{1}{N(\nu)^s}.
\]
We get:
\begin{cor}\label{fcteq}
The functional equation shows that 
$F(\bfr, \ffr,\epsilon,x,1- k)$ for $k=1,2,\ldots$ is non-zero
for $|\epsilon|=0$ and $k$ even or for $|\epsilon|=g$ and $k$ odd. In these
cases one has
\[
\zeta(\bfr, \ffr,x,1- k)=|d_F|^{-\frac{1}{2}}N(\ffr^{-1}\bfr)
\frac{((k-1)!)^g}{(2\pi i)^{kg}}L(\bfr, \ffr,\epsilon,x, k).
\]
\end{cor}
\subsection{The main theorem}
Here we formulate our main theorem. It computes the residue
map from (\ref{deg2}) in terms of 
the partial $L$-functions.

The transcendental description of the cusps gives
\[
H^0(\dS(\C),L)=\Ind_{B^1(\Z)}^{G(\Z/n\Z)}L
\]
and $H^0(\dS,L)$ is the subgroup of elements invariant under $D(\Z/n\Z)$.
Similarly, the $n$-torsion sections of $\Ah[n]$ over $S(\C)$
can be identified with functions from $G(\Z/n\Z)$ to $(\frac{1}{n}\Oh/\Oh)^2$,
which are equivariant with respect to the canonical 
$G^1(\Z):=\ker(G(\Z)\to \Z^*)$ action.
The action of $G(\Z/n\Z)$ on $S$ induces via  pull-back an action
on $\Ah[n](S(\C))$ and we have:
\[
\Ah[n](S(\C))=\Ind_{ G^1(\Z)}^{G(\Z/n\Z)}(\frac{1}{n}\Oh/\Oh)^2.
\]
The group $\Ah[n](S)$ consists again of the elements 
invariant under $D(\Z/nZ)$.
Let $D:=\Ah[n](S)$ and consider the formal linear combinations
\[
L[D]^0:=\bigl\{ \sum_{\sigma\in D}l_\sigma(\sigma):l_\sigma\in L\text{ and }
\sum_{\sigma\in D} l_\sigma=0\bigr\}.
\] 
The ${G(\Z/n\Z)}$ action on $D$ carries over to an action on 
$L[D]^0$.
For $\alpha\in L[D]^0$ and $k>0$ (or $\alpha\in L[D\setminus e(S)]^0$
and $k\geq 0$)
we construct in \ref{eisdefn} a class
\[
\Eis^k(\alpha)\in \Ext^{2g-1}_S(L,\Sym^k\Hh(g)),
\]
which depends on $\alpha$ in a functorial way. Thus, the resulting map 
\begin{equation}\label{resequiv}
L[D]^0\xrightarrow{\Eis^k}\Ext^{2g-1}_S(L,\Sym^k\Hh(g))\xrightarrow{\res}
\Ind_{B^1(\Z)}^{G(\Z/n\Z)}L
\end{equation}
is equivariant for the $G(\Z/n\Z)$ action. 
\begin{thm}\label{maintheorem} Let $L\supset \Q$ and
$\alpha=\sum_{\sigma\in D}l_\sigma(\sigma)$. 
Then $\res(\Eis^{m}(\alpha))$ is non-zero
only for $m\equiv 0(g)$ and for every $h\in G(\Z/n\Z)$ and $k>0$
\[
\res(\Eis^{gk}(\alpha))(h)=(-1)^{g-1}
\sum_{\sigma\in D}l_\sigma \zeta(\Oh,\Oh,p(h\sigma),-k).
\]
\end{thm}
To use the basis given by the coinvariants in $\Sym^{gk}\Oh\otimes L$
as we did in the proof of theorem \ref{degenercompu} is not natural.
A better description is as follows:
For each $h\in G(\Z/n\Z)$
choose an element $d_h\in D(\Z/n\Z)$ such that 
$\tilde{h}:=hd_h^{-1}\in G^1(\Z/n\Z)$. Then, as in (\ref{phdefn})
we have an ideal $\bfr_{\tilde{h}}$ and a projection
\[
\Oh^2\xrightarrow{p_\tildeh}\bfr_\tildeh.
\]
Now use the identification
$H^{g-1}(S_T, \Sym^{gk}\bfr_\tildeh\otimes L)\isom L$ at the cusp $h$. With
this basis the above result reads
\begin{cor} In this basis 
\[
\res(\Eis^{gk}(\alpha))(h)=(-1)^{g-1}N\bfr_\tildeh^{-k-1}
\sum_{\sigma\in D}l_\sigma \zeta(\bfr_\tildeh,\Oh,p_\tildeh(\sigma),-k).
\]
\end{cor}
The theorem and the corollary will be proved in section \ref{proof}.

\section{Polylogarithms}\label{polsection}
In this section we review the theory of the polylogarithm on abelian schemes.
Special emphasis is given the topological case, which will be important
in the proof of the main theorem.  
The elliptic polylogarithm was introduced by Beilinson and Levin \cite{BL}
and the generalization to higher dimensional families of abelian
varieties is due to Wildeshaus \cite{W}. The idea to interprete the construction
by Nori in terms of the topological polylogarithm is due to Beilinson and Nori
(unpublished).

The polylogarithm can be defined in any of the
categories $MHM, et, top$ for any abelian scheme $\pi:\Ah\to S$,
with unit section $e:S\to \Ah$ of constant relative dimension
$g$. If we work in $top$, it even suffices to assume that
$\pi:\Ah\to S$ is a family of topological tori (i.e., fiberwise isomorphic
to $(\R/\Z)^g$). 
For more details in the case of abelian schemes, see \cite{W} chapter III
part I, or \cite{L}. In the case of elliptic curves one can also consult
\cite{BL} or \cite{HK}.

\subsection{Construction of the polylog}
For simplicity we assume $L\supset \Q$ in this section and discuss the
necessary modifications for integral coefficients later. 
Define a lisse sheaf $\Log^{(1)}$ on $\Ah$, which is an extension
\[
0\to \Hh\to\Log^{(1)}\to L\to 0
\]
together with a splitting $s:e^*L\to e^*\Log^{(1)}$
in any of the three categories $MHM, et, top$ as follows: Consider the
exact sequence
\[
0\to \Ext^1_{S}(L,\Hh)\xrightarrow{\pi^*} \Ext^1_{\Ah}(L,\pi^*\Hh)\to 
\Hom_{S}(L,R^1\pi_*\pi^*\Hh)\to 0,
\]
which is split by $e^*$. Note that by the projection formula
$ R^1\pi_*\pi^*\Hh\isom R^1\pi_*L\otimes \Hh$ so that
\[
\Hom_S(L,R^1\pi_*\pi^*\Hh)\isom \Hom_S(\Hh, \Hh).
\]
Then $\Log^{(1)}$ is a sheaf representing
the unique extension class in $\Ext^1_{\Ah}(L,\pi^*\Hh)$,
which splits when pulled back to $S$ via $e^*$ and which maps to 
$\id \in \Hom_S(\Hh, \Hh)$. Define 
\[
\Log^{(k)}:=\Sym^k\Log^{(1)}.
\]
\begin{defn} The {\em logarithm sheaf} is the pro-sheaf
\[
\Log:=\Log_\Ah:=\prolim \Log^{(k)},
\]
where the transition maps are  induced by the map $\Log^{(1)}\to L$.
In particular, one has exact sequences
\[
0\to \Sym^k\Hh\to \Log^{(k)}\to \Log^{(k-1)}\to 0
\]
and a splitting induced by $s:e^*L\to e^*\Log^{(1)}$
\[
e^*\Log\isom \prod_{k\geq 0}\Sym^k\Hh.
\]
\end{defn}
Any isogeny $\phi:\Ah\to \Ah$ of degree invertible in $L$
induces an isomorphism $\Log\isom \phi^*\Log$, which is on the associated 
graded induced by $\Sym^k\phi:\Sym^k\Hh\to \Sym^k\Hh$.
For every torsion point $x\in \Ah(S)_{\tors}$ one gets an isomorphism
\begin{equation}\label{logtriv}
x^*\Log\isom e^*\Log\isom \prod_{k\geq 0}\Sym^k\Hh.
\end{equation}
The most important property of the sheaf $\Log$ is the vanishing of
its higher direct images except in the highest degree.
\begin{thm}[Wildeshaus, \cite{W}, cor. 4.4., p. 70]\label{logcohvanish} 
One has 
\[
R^i\pi_*\Log=0\mbox{ for }i\neq 2g
\]
and the augmentation $\Log\to L$ induces canonical isomorphisms
\[
R^{2g}\pi_*\Log\isom R^{2g}\pi_*L\isom L(-g).
\]
\end{thm}
For the construction of the polylogarithm 
one considers a non-empty disjoint union of torsion sections $i:D\subset \Ah$,
whose orders are invertible in $L$ 
(more generally, one can also consider $D$ {\'e}tale over $S$). 
Let 
\[
L[D]:=\bigoplus_{\sigma\in D}L
\]
and $L[D]^0\subset  L[D] $ the kernel of the augmentation map $L[D]\to L$.
Elements $\alpha\in L[D]$ are written as formal linear combinations
$\alpha=\sum_{\sigma\in D}l_{\sigma}(\sigma)$.
Similarly, define
\[
\Log[D]:=\bigoplus_{\sigma\in D}\sigma^*\Log
\]
and 
\[
\Log[D]^0:=\ker\left( \Log[D]\to L\right)
\]
to be the kernel of the composition of the 
sum of the augmentation maps $\Log[D]\to L[D]$ and the augmentation $L[D]\to L$.
\begin{cor} \label{eq:polDlocseq}
The localization sequence for $U:=\Ah\setminus D$ induces 
an isomorphism
\[
\Ext^{2g-1}_{U}(L[D]^0, \Log(g))\isom \Hom_S(L[D]^0,\Log[D]^0).
\]
\end{cor}
\bew
The vanishing result \ref{logcohvanish} implies that the localization
sequence is of the form 
\[
0\to \Ext^{2g-1}_{U}(L[D]^0, \Log(g))\to \Hom_S(L[D]^0,i^*\Log)\to 
\Hom_S(L[D]^0,L)\to 0.
\]
Inserting the definition of $\Log[D]^0$ gives the desired result.
\bewende

\begin{defn}\label{logdefn} 
The {\em polylogarithm} $\pol^D$ is the extension class
\[
\pol^D\in \Ext^{2g-1}_{U}(L[D]^0, \Log(g)),
\]
which maps to the canonical inclusion $L[D]^0\to \Log[D]$ 
under the isomorphism in \ref{eq:polDlocseq}.
In particular, for every $\alpha\in L[D]^0$ we get by pull-back an extension class
\[
\pol^D_\alpha \in \Ext^{2g-1}_{U}(L, \Log(g)).
\]
\end{defn}

\subsection{Integral version of the polylogarithm, the topological case}

In the topological and the \'etale situation it is possible 
to define the polylogarithm with integral coefficients. In this section
we treat the topological case and the \'etale case in the next section.
The construction presented here is a reinterpretation by Beilinson and
Levin (unpublished) of results of Nori and Sczech. 

We start by defining the logarithm sheaf 
for any (commutative) coefficient ring $L$, in particular for $L=\Z$.
In the topological situation, it is even possible to define
more generally 
the polylogarithm for any smooth family of real tori of constant dimension $g$,
which has a unit section.

Let
\[
\pi:\Th\to S
\]
be such a family, $e:S\to \Th$ the unit section 
and let $\Hh_L:= \underline{\Hom}_S(R^1\pi_*L,L)$ be
the local system of the homologies of the fibers with coefficients in $L$.
Let $\tilde{\Hh}_\R$ be the associated vector bundle of $\Hh_\R$.
Then $\Th\isom \Hh_\Z\backslash\tilde{\Hh}_\R$ and we denote by 
\[
\tilde{\pi}:\tilde{\Hh}_\R\to \Th
\]
the associated map.
Let 
\[
L[\Hh_\Z]:=e^*\tilde{\pi}_!L
\]
be the local system of group rings on $S$, which is stalk-wise
the group ring of the stalk of the local system $\Hh_\Z$ with coefficients
in $L$. The augmentation ideal of  $ L[\Hh_\Z]\to L$ is denoted by $\Ih$ 
and we define 
\[
L[[\Hh_\Z]]:=\prolim_r L[\Hh_\Z]/\Ih^r
\]
the completion along the augmentation ideal. 
Note that $\Ih^n/\Ih^{n+1}\isom \Sym^n\Hh_L$. If $L\supset \Q$,
one has even a ring isomorphism
\begin{equation}\label{logisom}
L [[\Hh_\Z]]\xrightarrow{\isom} \prod_{k\geq 0}\Sym^k\Hh_L,
\end{equation}
induced
by $h\mapsto \sum_{k\geq 0}h^{\otimes k}/k!$ for $h\in\Hh_\Z$. 
\begin{defn}
The {\em logarithm sheaf} $\Log$ is the local system on $\Th$
defined by
\[
\Log:=\tilde{\pi}_! L \otimes_{L[\Hh_\Z]}L[[\Hh_\Z]].
\]
As a local system of $L[[\Hh_\Z]]$-modules, $\Log$ is of rank $1$.
\end{defn}
Any isogeny $\phi:\Th\to\Th$ of order invertible in $L$ induces 
an isomorphism $\Log\isom\phi^*\Log$, which is induced by 
$\phi:\Hh_\Z\to \Hh_\Z$.  In particular,
if the order of a torsion section $x:S\to \Th$ is invertible in
$L$, one has an isomorphism
\[
x^*\Log\isom e^*\Log=L[[\Hh_\Z]].
\]
To complete the definition of the polylogarithm, one has to compute
the cohomology of $\Log$. As $L[[\Hh_\Z]]$ is a flat $L[\Hh_\Z]$-module
one gets
\[
R^i\pi_*\Log\isom R^i\pi_*\tilde{\pi}_!L\otimes_{L[\Hh_\Z]}L[[\Hh_\Z]]
\]
and because $\pi_*=\pi_!$ one has to consider
$ R^i(\pi\verk \tilde{\pi})_!L$. But the fibers of 
\[
\pi\verk\tilde{\pi}: \tilde{\Hh}_\R\to S
\]
are just $g$-dimensional vector spaces and the cohomology with compact support
lives only in degree $g$, where it is the dual of $\Lambda^{max}\Hh_L$. 
Hence, we have proved:
\begin{lemma} Denote by $\mu_\Th^\lor$
the $L$-dual of $\mu_\Th:=\Lambda^{max}\Hh_L$. 
Then the higher direct images of 
$\Log$ are given by
\[
R^i\pi_*\Log\isom \left\{\begin{array}{cc}\mu_\Th^\lor&\mbox{ if }i=g\\
0&\mbox{ else. }\end{array}\right.
\]
\end{lemma}
As in \ref{eq:polDlocseq} one obtains
\[
\Ext^{g}_{U}(L[D]^0, \Log\otimes\mu_\Th)\isom \Hom_S(L[D]^0,\Log[D]^0)
\]
and one defines the polylogarithm 
\[
\pol^D\in \Ext^{g-1}_{U}(L[D]^0, \Log\otimes\mu_\Th)
\]
in the same way. For $\alpha\in L[D]^0$ one has again
\[
\pol^D_\alpha\in \Ext^{g-1}_{U}(L, \Log\otimes\mu_\Th)=
H^{g-1}(U,\Log\otimes\mu_\Th).
\]
\subsection{Integral version of the polylogarithm, the \'etale case}

This section will not be used in the rest of the paper and can be omitted
by any reader not interested in the integral \'etale case. 

To define an integral \'etale polylogarithm, one has to modify the
definition of the logarithm sheaf as in the topological case. The 
situation we consider here is again an abelian scheme
\[
\pi:\Ah\to S
\]
of constant fiber dimension $g$ and unit section $e:S\to \Ah$.
Let $\ell$ be a prime number,
$L=\Z/\ell^k\Z$ and assume that $\ell$ is invertible in $\Oh_S$. 
Then the $\ell^r$-multiplication $[\ell^r]:\Ah\to \Ah$ is \'etale 
and the sheaves $[\ell^r]_!L$ form a projective system via the
trace maps
\[
[\ell^r]_!L\to [\ell^{r-1}]_!L.
\]
\begin{defn} The {\em logarithm sheaf} is the inverse limit
\[
\Log_L:=\prolim_r[\ell^r]_!L
\]
with respect to the above trace maps. The logarithm sheaf with
$\Z_\ell$-coefficients is defined by
\[
\Log_{\Z_\ell}:=\prolim_k\Log_{\Z/\ell^k\Z}.
\]
\end{defn}
Let $\Hh_\ell:=\prolim_r\Ah[\ell^r]$ 
be the Tate-module of $\Ah/S$. As $\ell$ is nilpotent in $L$, 
we get that $e^*\Log=L[[\Hh_\ell]]$ is the Iwasawa algebra of 
$\Hh_\ell$ with coefficients in $L$. 
Any isogeny $\phi: \Ah\to \Ah$ of degree prime to $\ell$
induces an isomorphism $[\ell^r]_!L\to \phi^*[\ell^r]_!L$, which induces
\[
\Log\isom \phi^*\Log.
\]
\begin{prop} Let $L=\Z/\ell^k\Z$ or $L=\Z_\ell$.
The higher direct images of 
$\Log$ are given by
\[
R^i\pi_*\Log\isom \left\{\begin{array}{cc}L(-g)&\mbox{ if }i=2g\\
0&\mbox{ else. }\end{array}\right.
\]
\end{prop}
\bew It suffices to consider the case $L=\Z/\ell^k\Z$. 
We will show that the transition maps 
$R^i\pi_*[\ell^r]_!L\to R^i\pi_*[\ell^{s}]_!L$
are zero for $i<2g$ and every $s$, if $r$ is sufficiently big.
By Poincar\'e duality we may consider the maps
\begin{equation}\label{dualtransition}
R^{2g-i}\pi_![\ell^s]_*L(g)\to R^{2g-i}\pi_![\ell^{r}]_*L(g).
\end{equation}
By base change we may assume that $S$ is the spectrum of an
algebraically closed field. Denote by $\Ah_s$ the variety $\Ah$ considered
as covering of $\Ah$ via $[\ell^s]$. Then
\[
R^1\pi_![\ell^s]_*L(g)=H^1(\Ah,[\ell^s]_*L(g))=\Hom(\pi_1(\Ah_s),L(g)).
\]
With this description we see that for every 
$f\in \Hom(\pi_1(\Ah_s),L(g))$ there is an $r$, such that the restriction
to $\pi_1(\Ah_r)$ is trivial. This shows that the map in (\ref{dualtransition})
is zero, if $r$ is sufficiently big and $i<2g$ as the cohomology in
degree $i$ is the $i$-th exterior power of the first cohomology. That 
(\ref{dualtransition}) is an isomorphism for $i=2g$ is clear.
\bewende

\subsection{Eisenstein classes}
The Eisenstein classes are specializations of the polylogarithm.
The situation is as follows. 
First let $\alpha\in L[\Ah[n]\setminus e(S)]^0$ 
and assume that $\Q\subset L$.
Then one can pull-back the class 
$\pol^{\Ah[n]\setminus e(S)}_\alpha\in \Ext^{g-1}_{U}(L, \Log(g))$
along $e$ and gets:
\[
e^*\pol^{\Ah[n]\setminus e(S)}_\alpha\in \Ext^{2g-1}_{S}(L, e^*\Log(g))=
\prod_{k\geq 0}\Ext^{2g-1}_{S}(L,\Sym^k\Hh(g)).
\]
The $k$-th component is the {\em Eisenstein class} $\Eis^k(\alpha)$.
For $k>0$, we can extend this definition to $\alpha\in L[\Ah[n]]^0$
with the following observation:
\begin{lemma} Let $\lambda\in\Oh$ and $[\lambda]:\Ah\to\Ah$ the associated
isogeny. Assume that the degree of $[\lambda]$ is prime to $n$. Then
$[\lambda]$ induces via pull-back an isomorphism
\[
\lambda^*:L[\Ah[n]]^0\to L[\Ah[n]]^0
\]
and for $k>0$
\[
\Eis^k(\lambda^*(\alpha))=\lambda^k\Eis^k(\alpha).
\]
Here the $\lambda^k\Eis^k(\alpha)$ uses the  $\Oh$-module structure
on $\Ext^{2g-1}_{S}(L,\Sym^k\Hh(g))$.
\end{lemma}
\bew It is clear that $\lambda^*$ is an isomorphism.
By definition $\pol^{\Ah[n]\setminus e(S)}_\alpha$ is functorial
with respect to isogenies and one only has to remark that 
$\Log\isom [\lambda]^*\Log$ is on the associated graded given by
$\Sym^\bullet[\lambda]:\Sym^\bullet\Hh\to \Sym^\bullet\Hh$.
\bewende
Let now $\alpha\in L[\Ah[n]]^0$, then for $\lambda\neq 1,0$
$\alpha-\lambda^*\alpha L[\Ah[n]\setminus e(S)]^0$ and we define for $k>0$
\begin{equation}\label{eisen}
\Eis^k(\alpha):=(1-\lambda^k)^{-1}\Eis^k(\alpha-\lambda^*\alpha).
\end{equation}
It is a straightforward computation, that this definition does
not depend on the choosen $\lambda$.
\begin{defn}\label{eisdefn}
For any $\alpha\in L[\Ah[n]\setminus e(S)]^0$, 
define the $k$-th {\em Eisenstein class} associated to $\alpha$,
\[
\Eis^k(\alpha)\in \Ext^{2g-1}_{S}(L,\Sym^k\Hh(g)),
\]
to be the $k$-th component of $e^*\pol^{D}_\alpha$. For $k>0$ and
$\alpha\in L[\Ah[n]]^0$ define $\Eis^k(\alpha)$ by the formula in
(\ref{eisen}).
\end{defn}
Note that by the functoriality of the polylogarithm the map
\begin{equation}\label{eisenequiv}
L[\Ah[n]\setminus e(S)]^0\xrightarrow{\Eis^k}\Ext^{2g-1}_{S}(L,\Sym^k\Hh(g))
\end{equation}
is equivariant for the $G(\Z/n\Z)$ action on both sides.

These Eisenstein classes should be considered
as analogs of Harder's Eisenstein classes (but observe that
we have only classes in cohomological degree $2g-1$). The advantage of the
above classes is that they are defined by a universal condition, which
makes a lot of their properties  easy to verify.

\section{Proof of the main theorem}\label{proof}
In this section we assume that $\Q\subset L$.

The proof of the main theorem will be in several steps.
First we reduce to the case of local systems for the usual topology.
The second step consists of a trick already used in \cite{HK}: instead
of working with the Eisenstein classes directly, we work with the polylogarithm
itself. The reason is that the polylog is characterized by a universal 
property and has a very good functorial behavior.
The third step reviews  the computations of Nori in \cite{Nori}.
In the fourth step we compute the integral over $S^1_T$ and
the fifth step gives the final result.
\subsection{1. Step: Reduction to the classical topology}
We distinguish the $MHM$ and the \'etale case. 
In the $MHM$ case, the target of the residue map from (\ref{deg2})
\begin{equation}
\res:\Ext^{2g-1}_S(L, \Sym^k\Hh(g))\to\Hom_{\dS}(L,L).
\end{equation}
is purely topological and does not depend on the Hodge structure. 
More precisely, 
the canonical map ``forget the Hodge structure'' denoted by $\rat$
induces an isomorphism
\[
\rat:\Hom_{MHM,\dS}(L,L)\isom \Hom_{top,\dS}(L,L).
\]
By \cite{Sa} thm. 2.1 we have  a commutative diagram
\begin{equation}
\begin{CD} \Ext^{2g-1}_{MHM,S}(L, \Sym^k\Hh(g))@>\res >>\Hom_{MHM,\dS}(L,L)\\
@VV\rat V@VV \rat V\\
\Ext^{2g-1}_{top,S}(L, \Sym^k\Hh(g))@>\res >>\Hom_{top,\dS}(L,L).
\end{CD}
\end{equation}
This reduces the computation of the residue map for $MHM$ to the case
of local systems in the classical topology.

In the \'etale case one has an injection
\[
\Hom_{et,\dS}(L,L)\hookrightarrow \Hom_{et,\dS\times\Qbar}(L,L)\isom
\Hom_{top,\dS(\C)}(L,L).
\]
and a commutative diagram
\begin{equation}
\begin{CD} \Ext^{2g-1}_{et,S}(L, \Sym^k\Hh(g))@>\res >>\Hom_{et,\dS}(L,L)\\
@VVV@VVV\\
\Ext^{2g-1}_{top,S(\C)}(L, \Sym^k\Hh(g))@>\res >>\Hom_{top,\dS(\C)}(L,L).
\end{CD}
\end{equation}
Again, this reduces the residue computation to the classical topology.

\subsection{2. Step: Topological degeneration} 
In this section we reduce the computation of $\res\circ\Eis^k$ to 
a computation of the polylog on $\Th_\Mh$.

We are now in the topological situation and use again the 
notations $\dS$ and $S$ instead of $\dS(\C)$ and $S(\C)$.

Recall from (\ref{resequiv}) that $\res\circ\Eis^k$ is $G(\Z/n\Z)$
equivariant. In particular,
\[
\res(\Eis^k(\alpha))(h)=\res(\Eis^k(h\alpha))(\id),
\]
where $h\alpha$ denotes the action of $h$ on $\alpha$. To compute
the residue it suffices to consider the residue at $\id$.

Recall from (\ref{degenerationdiagr}) that
we have a commutative diagram of fibrations
\begin{equation}\label{fibrediagr}
\begin{CD}
\Ah(\C)@>p>>\Th_\Mh\\
@V\pi VV@VV\pi_\Mh V\\
S^1_B @>q>> S^1_T.
\end{CD}
\end{equation}
The map $p:\Hh\to \Mh$ induces $\Log_\Ah\to p^*\Log_{\Mh}$.
Let $D=\Ah[n]$ and  $U:=\Ah\setminus D$ be the complement.
Let $p(D)=\Th_\Mh[n]$ be the image of $D$ in $\Th_\Mh$ and
$V:=\Th_\Mh\setminus p(D)$ be its complement in $\Th_\Mh$.
Then $p$ induces a map 
\[
p:U\setminus p^{-1}(p(D))\to V.
\]
We define a trace map
\begin{equation}
  \label{eq:poltrace}
  p_*:\Ext^{2g-1}_{U}(L, \Log_\Ah\otimes\mu_\Ah)\to 
\Ext^{g-1}_{V}(L, \Log_{\Mh}\otimes\mu_{\Th_\Mh})
\end{equation}
as the composition of the restriction to $U\setminus p^{-1}(p(D))$
\[
\Ext^{2g-1}_{U}(L, \Log_\Ah\otimes\mu_\Ah)\to 
\Ext^{2g-1}_{U\setminus p^{-1}(p(D))}(L, \Log_\Ah\otimes\mu_\Ah)
\]
with the adjunction map
\[
\Ext^{2g-1}_{U\setminus p^{-1}(p(D))}(L, \Log_\Ah\otimes\mu_\Ah)\to 
\Ext^{g-1}_{V}(L,R^gp_*p^*\Log_\Mh\otimes\mu_\Ah).
\]
As $\mu_\Ah\isom\mu_{\Th_\Nh}\otimes\mu_{\Th_\Mh}$, 
the projection formula gives 
\[
R^gp_*p^*\Log_\Mh\isom \Log_\Mh\otimes\mu_{\Th_\Mh}.
\] 
The composition of these maps gives the desired $p_*$ in (\ref{eq:poltrace}).
The crucial fact is that the polylogarithm behaves well under
this trace map.
\begin{prop}\label{poldeg} With the notations above, 
let $\alpha\in L[D]^0$ and 
$\pol^D_{\Ah,\alpha}\in \Ext^{2g-1}_{U}(L, \Log_\Ah)$ be the associated 
polylogarithm. Denote by $p(\alpha)$ the image of $\alpha$ under the map
\begin{align*}
p:L[D]^0&\to L[p(D)]^0
\end{align*}
induced by $p:\Ah(\C)\to \Th_\Mh$. 
Then 
\[
p_*\pol^D_{\Ah,\alpha}=\pol^{p(D)}_{\Th_\Mh, p(\alpha)}.
\]
\end{prop}
\bew This is a quite formal consequence of the definition and the fact that
the residue map commutes with the trace map.
We use cohomological notation, then one has a commutative diagram
\[\begin{CD} 
H^{2g-1}(U,\Log_\Ah\otimes\mu_\Ah)@>>> H_D^{2g}(\Ah,\Log_\Ah\otimes\mu_\Ah)\\
@VVV@VVV\\
H^{2g-1}(U\setminus p^{-1}(p(D)),\Log_\Ah\otimes\mu_\Ah)@>>> 
H_{p^{-1}(p(D))}^{2g}(\Ah,\Log_\Ah\otimes\mu_\Ah)\\
@VVp_*V@VVp_*V\\
H^{g-1}(V,\Log_\Mh\otimes\mu_{\Th_\Mh})@>>> 
H_{p(D)}^{g}(\Th_\Mh,\Log_\Mh\otimes\mu_{\Th_\Mh}).
\end{CD}
\]
We can identify
\[
H_D^{2g}(\Ah,\Log_\Ah\otimes\mu_\Ah)\isom \bigoplus_{\sigma\in D} 
\sigma^*\Log_\Ah
\]
and 
\[
H_{p(D)}^{g}(\Th_\Mh,\Log_\Mh\otimes\mu_{\Th_\Mh})\isom
\bigoplus_{\sigma\in p(D)}\sigma^*\Log_{\Th_\Mh}.
\]
With this identification the composition of the vertical arrows
on the right is induced by $\Log_\Ah\to p^*\Log_{\Th_\Mh}$. The polylog
$\pol^D_{\Ah,\alpha}$ belongs to the section 
$\alpha\in L[D]^0\subset\bigoplus_{\sigma\in D} \sigma^*\Log_\Ah$. 
This maps to 
$p(\alpha)\in L[p(D)]^0\subset \bigoplus_{\sigma\in p(D)}
\sigma^*\Log_{\Th_\Mh}$.
Thus $\pol^D_{\Ah,\alpha}$ is 
mapped under
$p_*$ to $\pol^{p(D)}_{\Th_\Mh, p(\alpha)}$.
\bewende

We want to prove the same sort of result for the Eisenstein classes
themselves. To formulate it properly, we need:
\begin{lemma}\label{GammaNcoinv} Let $q:S^1_B\to S^1_T$ 
be the fibration from (\ref{fibrediagr}).
Then 
\[
R^gq_*\Sym^k\Hh\isom \Sym^k\Mh\otimes\mu_{\Th_\Nh}^{\lor}.
\]
\end{lemma}
\bew Recall the exact sequence
\[
0\to \Nh\to\Hh\to\Mh\to 0
\]
from (\ref{locexseq}). By definition of $N(\Z)$, the coinvariants
of $\Sym^k\Hh$ for $N(\Z)$ are exactly $\Sym^k\Mh$.
The lemma follows, as $R^gq_*$ corresponds by definition of the
fibering exactly to the coinvariants under $N(\Z)$.
\bewende 
Define a trace map 
\[
q_*:\Ext^{2g-1}_{S_B}(L,\Sym^k\Hh\otimes\mu_\Ah)\to
\Ext^{g-1}_{S_T}(L,\Sym^k\Mh\otimes\mu_{\Th_\Mh})
\]
by adjunction for $q$, the 
isomorphism $R^gq_*\Sym^k\Hh\isom\Sym^k\Mh\otimes\mu_{\Th_\Nh}^{\lor}$
from lemma \ref{GammaNcoinv}
and the isomorphism $\mu_\Ah\isom\mu_{\Th_\Nh}\otimes\mu_{\Th_\Mh}$.
The behaviour of $\Eis^k(\alpha)$ under $q_*$ is given by:
\begin{thm}\label{eisendeg}
Let $k>0$ and $\alpha\in L[D]^0$. Then
\[
q_*(\Eis^k_\Ah(\alpha))=\Eis^k_{\Th_\Mh}(p(\alpha)),
\]
where $p:L[D]^0\to L[p(D)]^0$ is the map from \ref{poldeg}.
\end{thm}
\bew Consider the following diagram in the derived 
category:
\begin{equation}\label{commdiagr}
\begin{CD}Rp_*\Log_\Ah \otimes\mu_\Ah @>>> Rp_*e_*e^*\Log_\Ah \otimes\mu_\Ah\\
@VVV@|\\
Rp_*p^*\Log_{\Th_\Mh}\otimes\mu_\Ah @. e'_*Rq_*e^*\Log_\Ah \otimes\mu_\Ah\\
@VVV@VVV\\
\Log_{\Th_\Mh}\otimes\mu_{\Th_\Mh}[-g]@>>> 
e'_*e'^*\Log_{\Th_\Mh}\otimes\mu_{\Th_\Mh}[-g]
\end{CD}
\end{equation}
We will show that this diagram is commutative and thereby explain all the
maps. First consider the commutative diagram
\[
\begin{CD}Rp_*\Log_\Ah \otimes\mu_\Ah @>>> Rp_*e_*e^*\Log_\Ah \otimes\mu_\Ah\\
@VVV@VVV\\
Rp_*p^*\Log_{\Th_\Mh}\otimes\mu_\Ah @>>>
Rp_*e_*e^*p^*\Log_{\Th_\Mh}\otimes\mu_\Ah,
\end{CD}
\]
where the horizontal arrows are induced from adjunction $\id\to e_*e^*$
and the vertical arrows from $ \Log_\Ah \to p^*\Log_{\Th_\Mh}$.
One has $p\circ e= e'\circ q$ and hence
\[
Rp_*e_*e^*p^*\Log_{\Th_\Mh}\otimes\mu_\Ah\isom 
e'_*Rq_*q^*e'^*\Log_{\Th_\Mh}\otimes\mu_\Ah.
\]
The projection formula gives
\[
e'_*Rq_*q^*e'^*\Log_{\Th_\Mh}\otimes\mu_\Ah\isom
e'_*e'^*\Log_{\Th_\Mh}\otimes\mu_\Ah\otimes Rq_*L.
\]
Projection to the highest cohomology gives a commutative diagram
\[
\begin{CD}
Rp_*p^*\Log_{\Th_\Mh}\otimes\mu_\Ah @>>>
e'_*e'^*\Log_{\Th_\Mh}\otimes\mu_\Ah\otimes Rq_*L\\
@VVV@VVV\\
\Log_{\Th_\Mh}\otimes\mu_\Ah\otimes \mu_{\Th_\Nh}^\lor
@>>>
e'_*e'^*\Log_{\Th_\Mh}\otimes\mu_\Ah\otimes \mu_{\Th_\Nh}^\lor,
\end{CD}
\]
where the horizontal maps are adjunction maps $\id\to e'_*e'^*$.
Finally we use $\mu_\Ah\otimes \mu_{\Th_\Nh}^\lor\isom \mu_{\Th_\Mh}$ to
obtain the commutative diagram (\ref{commdiagr}).
Applying $\Ext^{2g-1}_{V}(L,-)$ to this diagram, where 
$V:=\Th_\Mh\setminus p(D)$ we get
\[
\begin{CD}
\Ext^{2g-1}_{U}(L,\Log_\Ah\otimes\mu_\Ah)@>>> 
\Ext^{2g-1}_{S}(L,e^*\Log\otimes\mu_\Ah)\\
@Vp_*VV@VVq_*V\\
\Ext^{g-1}_{V}(L,\Log_{\Th_\Mh}\otimes\mu_{\Th_\Mh})@>>> 
\Ext^{g-1}_{S}(L,e'^*\Log_{\Th_\Mh}\otimes\mu_{\Th_\Mh}).
\end{CD}
\]
Now, as $k>0$, we may assume that $\alpha\in L[D\setminus e(S)]^0$ and
$p(\alpha)\in L[p(D)\setminus e'(S)]^0$. The result follows then from 
proposition \ref{poldeg}.
\bewende
In a similar (but simpler) way one shows:
\begin{thm}\label{isogenyinvariance} 
Let $\phi:\Th_\Mh\to \Th_{\Mh'}$ be an isogeny of tori,
then  $\phi$ induces a morphism $\phi_*:e^*\Log_\Mh\to e^*\Log_{\Mh'}$ and
\[
\phi_*\Eis^k_{\Th_\Mh}(\alpha)=\Eis^k_{\Th_{\Mh'}}(\phi(\alpha)).
\]
\end{thm}
\subsection{3. Step: Explicit description of the polylog}
In this section we follow Nori \cite{Nori} to describe the polylog
$\pol^{\Th_\Mh[n]}_\beta$ for any $\beta\in L[\Th_\Mh[n]\setminus 0]^0$
explicitly. The presentation is also influenced by unpublished notes
of Beilinson and Levin.

In fact it is useful for the connection with $L$-functions to consider 
a more general situation and to allow arbitrary fractional
ideals $\afr$ instead just $\Oh$.
 
We assume $L=\C$. The geometric situation is this:
Recall that $T^1(\Z)=\Oh^*$ and let $\afr\subset F$ be a fractional
ideal with the usual $T^1(\Z)$-action. We can form as 
usual the semi direct product 
\[
\afr\rtimes T^1(\Z),
\]
where the multiplication is given by the formula 
$(v,t)(v',t')=(v+tv',tt')$.
Similarly, we can form $\afr\otimes\R\rtimes T^1(\R)$ and we define
\[
\Th_\afr:=\afr\rtimes T(\Z)\backslash 
\left(\afr\otimes\R\rtimes T^1(\R)\right)/K^T_\infty.
\]
We have 
\[
\pi_\afr:\Th_\afr\to S^1_T
\]
and we consider the polylog for this real torus bundle of relative dimension
$g$. The case $\Th_\Mh$ is the one where 
$\bigl(\begin{smallmatrix}a&0\\0&d\end{smallmatrix}\bigr)\in T^1(\Z)$  acts 
via $d\in \Oh^*$ on $\Oh$. 
Let us describe the logarithm sheaf $\Log_{\Th_\afr}$ in this setting.
As the coefficients are $L=\C$, we can use the isomorphism from (\ref{logisom})
\begin{align}
\C [[\afr]]&\xrightarrow{\isom}\prod_{k\geq 0}\Sym^k\afr_\C=:\Uhat(\afr) \\
\nonumber v&\mapsto \exp(v):=\sum_{k=0}^\infty \frac{v^{\otimes k}}{k!}
\end{align}
The action of 
$(0,t)\in\afr\rtimes T^1(\Z)$ on $\Uhat(\afr)$ is induced by
the action of $T^1(\Z)$ on $\afr$. 
The action of 
\[
(v,\id)\in\afr\rtimes T^1(\Z)
\]
on $\Uhat(\afr)$ is given by multiplication
with $\exp(v)$.
The logarithm sheaf $\Log_{\Th_\afr}$ is just the local system
defined by the quotient
\[
\afr\rtimes T^1(\Z)\backslash 
\left(\afr\otimes\R\rtimes T^1(\R)\times\Uhat(\afr) \right)/K^T_\infty.
\]
A $\Ch^\infty$-section $f$ of $\Log_{\Th_\afr}$
is a function $f:\afr\otimes\R\rtimes T^1(\R)\to\Uhat(\afr)$,
which has the equivariance property 
\[
f((v,t)(v',t'))=(v,t)^{-1}f(v',t').
\]
In a similar way, we can describe $\Log_{\Th_\afr}$-valued
currents. The global $\Ch^\infty$-section
\[
\exp(-v):(v,t)\mapsto \sum_{k=0}^\infty\frac{(-v)^{\otimes k}}{k!},
\]
with $(v,t)\in \afr\otimes \R\rtimes T^1(\R)$ defines a trivialization of
$\Log_{\Th_\afr}$ as $\Ch^\infty$-bundle. Every current  $\mu(v,t)$ 
with values in  $\Log_{\Th_\afr}$ can then be written in the form
\[
\mu(v,t)=\nu(v,t)\exp(-v),
\]
where $\nu(v,t)$ is now a current with values in the constant bundle
$\Uhat(\afr)$. In particular, $ \nu(v,t)$ is invariant under the action
of $\afr\subset \afr\rtimes T^1(\Z)$.
\begin{lemma} Let ${\bf v}:\afr\otimes\R\to \Uhat(\afr)$ be the
canonical inclusion given by 
$\afr\otimes\R\subset\Sym^1 \afr\otimes\C$, then the 
canonical connection $\nabla$ on $\Log_{\Th_\afr}$ acts on $\nu$ by
\[
\nabla \nu=(d-d{\bf v})\nu.
\]
\end{lemma}
\bew Straightforward computation.
\bewende

Following Nori \cite{Nori} we describe the polylog 
as a $\Log_{\Th_\afr}$-valued current $\mu(v,t)$ on 
$\Th_\afr$, such that
\begin{equation}\label{polproperty}
\nabla\mu(v,t)=\delta_\beta,
\end{equation}
where 
\[
\delta_\beta:=\sum_{\sigma\in D}l_\sigma\delta_\sigma
\]
and $\delta_\sigma$ are the currents defined by integration over the
cycles on $\Th_\afr$ given by the section $\sigma$.
If we write as above
\[
\mu(v,t)=\nu(v,t)\exp(-v)
\]
we get the equivalent condition
\begin{equation}\label{trivpolproperty}
(d-d{\bf v})\nu(v,t)=\delta_\beta.
\end{equation}
As $\nu(v,t)$ is invariant under the $\afr$-action,
we can develop $\nu(v,t)$ into a Fourier series
\begin{equation}\label{fourierseries}
\nu(v,t)=\sum_{\rho\in \afr^\lor}\nu_{\rho}(t)e^{2\pi i \rho(v)}.
\end{equation}
The property (\ref{trivpolproperty}) reads for the Fourier 
coefficients $\nu_{\rho}(t)$:
\begin{equation}
(d+2\pi id\rho-d\underline{v})\nu_\rho(t)=
(e^{-2\pi i \rho(\beta)})\vol,
\end{equation}
where $\vol$ is the unique constant coefficient $g$-form on $\afr\otimes\R$, 
such that the integral $\int_{\afr\otimes\R/\afr}\vol=1$ and
\[
e^{-2\pi i \rho(\beta)}:= \sum_\sigma l_\sigma e^{-2\pi i \rho(\sigma)}.
\]
We do not explain
in detail the method of Nori to solve this equation, we just give the
result. This suffices, because the cohomology class of
the polylogarithm is uniquely determined by the equation (\ref{polproperty})
and we just need to give a solution for it. 

Fix a positive definite
quadratic form $q$ on $\afr\otimes\R$, viewed
as an isomorphism
\[
q:(\afr\otimes\R)^\lor\isom \afr\otimes\R.
\]
Define a left action of $t\in T^1(\R)$ by $q_t(v,w):=q(t^{-1}v,t^{-1}w)$.
Consider $\rho$ as element in
$(\afr\otimes\R)^\lor$. Then $q_t(\rho)$ can be considered as a vector field 
and we denote by $\iota_\rho$ the contraction
with this vector field $q_t(\rho)$. 
We may also consider $q_t(\rho)$ as
element in $\Uhat(\afr)$ and denote this by ${\bf q}_t(\rho)$.
\begin{thm}[Nori]
With the notations above, one has for $0\neq \rho$
\[
\nu_\rho(t)=\sum_{m=0}^{g-1}
\frac{(-1)^m(e^{-2\pi i\rho(\beta)})}
{(2\pi i\rho(q_t(\rho))-{\bf q}_t(\rho))^{m+1}}
\iota_\rho(d\circ\iota_\rho)^m
\vol
\]
and 
\[
\nu_0(t)=0
\]
\end{thm}
\bew Write $\Phi_\rho$ for the operator multiplication by
$2\pi id\rho-d{\bf v}$ and $\Psi_\rho:=d+\Phi_\rho$. One checks that
$\Psi_\rho\circ \Psi_\rho=0=\iota_\rho\circ\iota_\rho$ and
that $\Psi_\rho\circ\iota_\rho+\iota_\rho\circ\Psi_\rho$ is an isomorphism.
Indeed $\Phi_\rho\circ\iota_\rho +\iota_\rho\circ\Phi_\rho$ is
multiplication by $2\pi i\rho(q_t(\rho))-{\bf q}_t(\rho)$ 
and $\Lh_\rho:=d\circ\iota_\rho+\iota_\rho\circ $ is the Lie 
derivative with respect
to the vector field $q_t(\rho)$. The formula in the theorem is
just 
\[
\iota_\rho\circ(\Psi_\rho\circ\iota_\rho+\iota_\rho\circ\Psi_\rho)^{-1}(e^{-2\pi i\rho(\beta)})\vol
\]
and to check that 
\[
\Psi_\rho\circ\iota_\rho\circ(\Psi_\rho\circ\iota_\rho+\iota_\rho\circ\Psi_\rho)^{-1}=\id
\]
note that $\iota_\rho\circ\Psi_\rho$ commutes with
$(\Psi_\rho\circ\iota_\rho+\iota_\rho\circ\Psi_\rho)^{-1}$ and
$\iota_\rho\circ\Psi_\rho(e^{-2\pi i\rho(\beta)})\vol=0$.
\bewende

\begin{cor}\label{polcurrent} 
The polylogarithm $\pol^{\Th_\afr[n]}_\beta$ 
is given in the topological realization by the current
\[
\mu(v,t)=\nu(v,t)\exp(-v)
\]
where $\nu(v,t)$ is the current given by
\[
\sum_{m=0}^{g-1}\sum_{k=0}^\infty{k+m\choose k}
\sum_{\rho\in \afr^\lor\setminus 0}
\frac{(-1)^me^{2\pi i\rho(v-\beta)}}
{(2\pi i\rho(q_t(\rho)))^{k+m+1}}{\bf q}_t(\rho)^{\otimes k}
\iota_\rho(d\circ\iota_\rho)^m\vol.
\]
\end{cor}
\bew This follows from the formula 
$\frac{1}{(A-B)^{m+1}}=\sum_{k=0}^\infty\frac{B^{k}}{A^{k+m+1}}
{k+m\choose k}$.
\bewende
The Eisenstein classes are obtained by pull-back of this current along the
zero section $e$. As for $k>0$ the series over the $\rho$ converges 
absolutely, this is defined and only terms with $m=g-1$ survive.
We get the following formula for the Eisenstein classes.
\begin{cor} Let $\beta\in \C[\Th_\afr[n]\setminus 0]^0$ 
and $k>0$, then the topological Eisenstein class is given
by 
\[
\Eis^k(\beta)=\frac{(k+g-1)!}{k!}\sum_{\rho\in \afr^\lor\setminus 0}
\frac{(-1)^{g-1}e^{-2\pi i\rho(\beta)}}
{(2\pi i\rho(q_t(\rho)))^{k+g}}{\bf q}_t(\rho)^{\otimes k}
q_t(\rho)^*\iota_{\Eh}\vol.
\]
Here, we have written $\Eh$ 
for the Euler vector field and $q_t(\rho)$ is considered
as a function $q_t(\rho):S_T\to \afr\otimes \R$, which maps $t$ to
the vector  $q_t(\rho)$.
\end{cor}
\bew From \ref{polcurrent} we have to compute 
\[
e^*\iota_\rho(d\circ\iota_\rho)^m\vol.
\]
For this remark that  the Lie derivative 
$\Lh_\rho= d\circ\iota_\rho + \iota_\rho\circ d$ with respect to
the vector field $q_t(\rho)$ acts in the same way on
$\vol$ as $d\circ\iota_\rho$. One sees immediately that 
$e^*\iota_\rho(d\circ\iota_\rho)^m\vol=0$,
if $m<g-1$ and 
a direct computation in coordinates gives 
that $\iota_\rho(\Lh_\rho)^{g-1}\vol=(g-1)!q_t(\rho)^*\iota_{\Eh}\vol$.
\bewende

\subsection{4. Step: Computation of the integral}
To finish the proof of theorem \ref{maintheorem} we have to compute
$u_*\Eis^k(\beta)$, where $u:S^1_T\to pt$ is the structure map.
As we need only to compute the corresponding integral for the
component of $S^1_T$ corresponding to $\id$, we let
$\Gamma_T\subset T^1(\Z)$ be the stabilizer of $\id \in T(\Z/n\Z)$ and
consider 
\[
u_{\id}:\Gamma_T\backslash \bigl(T^1(\R)/K^{T^1}_\infty\bigr) \to pt.
\]
To compute the integral, we introduce
coordinates on $T^1(\R)\isom (F\otimes \R)^*$ and on the torus $\Th_\afr$.
We identify $F\otimes \R\isom \prod_{\tau:F\to \R}\R$ and denote by
$e_1,\ldots,e_g$ the standard basis on the right hand side and by
$x_1,\ldots,x_g$ the dual basis. For any element
$u=\sum u_ie_i$ or $u=\sum u_ix_i$ we write $Nu:=u_1\cdots u_g$.
Let $q$ be the
quadratic form given by $\sum x_i^2$. We identify the orbit of $q$ under
$T^1(\R)$ with $(F\otimes \R)^*_+$ by mapping
\begin{align}
(F\otimes \R)^*&\to T^1(\R)q\\
\nonumber t&\mapsto q_{t}.
\end{align}
This map factors over $(F\otimes \R)^*_+$ and
the map is compatible with the $T^1(\Z)$ action on
both sides. We let $t_1,\ldots, t_g$ be coordinates on
$(F\otimes \R)^1$ so that $t_1^2,\ldots,t_g^2$ are coordinates on 
$(F\otimes \R)^1_+$.
If we write $\rho=\sum\rho_ix_i$ and $t_i:=x_i(t)$, then
\[
\rho(q_{t}(\rho))= \sum t_i^2\rho_i^2
\]
and ${\bf q}_{t}(\rho)$ has coordinates $t_i^2\rho_i$. More precisely,
if we let ${\bf e}_1,\ldots,{\bf e}_g$ be the basis
$e_1,\ldots,e_g$ considered as elements 
of $\Uhat(\afr)$, which identifies $\Uhat(\afr)$ with the
power series ring $\C[[{\bf e}_1,\ldots,{\bf e}_g]]$, then
${\bf q}_{t}(\rho)=\sum t_i^2\rho_i{\bf e}_i$.
The volume form is given by 
\[
\vol=|d_F|^{-1/2}N\afr^{-1}dx_1\land\ldots\land dx_g
\]
and we can write the Euler vector field as $\Eh=\sum x_i\partial_{x_i}$.
One gets (observe that $Nt=1$)
\[
q_{t}(\rho)^*
\iota_\Eh\vol=|d_F|^{-1/2}2^{g-1}N(\rho)N\afr^{-1}
\sum_{k=1}^g(-1)^{k-1}t_kdt_1\land\ldots
\widehat{dt_k}\ldots\land dt_g.
\]
Explicitly, the Eisenstein class is given as a current on $T^1(\R)$ by
\begin{multline}
\Eis^k(\beta)(t)=\\
\frac{(k+g-1)!}{k!}
\sum_{\rho\in\afr^\lor\setminus 0}
\frac{(-1)^{g-1}e^{-2\pi i\rho(\beta)}(\sum t_i^2\rho_i{\bf e}_i)^{\otimes k}}
{(2\pi i\sum\rho_i^2t_i^2)^{k+g}}
q_{t}(\rho)^*
\iota_\Eh\vol
\end{multline}
Define an isomorphism $(\R\otimes F)^1\times \R^*\isom(\R\otimes F)^*$
by mapping $(t,r)\mapsto y:=rt$.
Then we get:
\begin{equation}\label{volrelation}
\frac{dy_1}{y_1}\land\ldots\land\frac{dy_g}{y_g}=
\frac{dr}{r}\land
\sum_{k=1}^g(-1)^{k-1}t_kdt_1\land\ldots\widehat{dt_k}
\ldots\land dt_g.
\end{equation}
We use this decomposition to write $\Eis^k(\beta)(t)$ as a Mellin transform:
\begin{multline}
\Eis^k(\beta)(t)=\\
\sum_{\rho\in \afr^\lor\setminus 0}
{(-1)^{g-1}e^{-2\pi i\rho(\beta)}}
\int_{\R_{>0}}\displaylimits 
{e^{-u(2\pi i\sum\rho_i^2t_i^2)}}
\frac{(\sum t_i^2\rho_i{\bf e}_i)^{\otimes k}}{k!}
u^{k+g}\frac{du}{u}
\land q_{t}(\rho)^*
\iota_\Eh\vol.
\end{multline}
Substitute $u=r^2=N(y)^{2/g}$ and use (\ref{volrelation}) to get
\begin{multline}
\Eis^k(\beta)(t)=\\
\sum_{\rho\in \afr^\lor\setminus 0}
\frac{(-1)^{g-1}2^{g}e^{-2\pi i\rho(\beta)}N(\rho)}{|d_F|^{1/2}N\afr}
\int_{\R_{>0}}\displaylimits 
e^{-2\pi i\sum\rho_i^2y_i^2}
\frac{(\sum y_i^2\rho_i{\bf e}_i)^{\otimes k}}{k!}
N(y){dy_1}\land\ldots\land{dy_g}.
\end{multline}
The application of $u_{\id,*}$ amounts to integration over
\[
\Gamma_T\backslash \bigl(T^1(\R)/K^{T^1}_\infty\bigr)\isom 
\Oh^*_{(n)}\backslash (F\otimes \R)^1_+,
\]
where $\Oh^*_{(n)}$ are the totally positive units, which are
congruent to $1$ modulo the ideal generated by $(n)$. 
This gives with the usual trick 
\begin{multline}
u_{\id,*}\Eis^k(\beta)=\\
\sum_{\rho\in \Oh^*_{(n)}\backslash\bigl(\afr^\lor\setminus 0\bigr)}
\frac{(-1)^{g-1}2^ge^{-2\pi i\rho(\beta)}N(\rho)}{|d_F|^{1/2}N\afr}
\int_{(F\otimes\R)^*_+}\displaylimits 
e^{-2\pi i\sum\rho_i^2y_i^2}
\frac{(\sum y_i^2\rho_i{\bf e}_i)^{\otimes k}}{k!}
N(y){dy_1}\land\ldots\land{dy_g}.
\end{multline}
The integral is a product of integrals for $j=1,\ldots, g$:
\[
\int_{\R_{>0}}\displaylimits 
e^{-2\pi i\rho_j^2y_j^2}\rho_j^k\frac{{\bf e}_j^{\otimes k}}{k!}y_j^{2k+2}
\frac{dy_j}{y_j}=
\frac{{\bf e}_j^{\otimes k}}
{2\rho_j(2\pi i\rho_j)^{k+1}}.
\]
We now consider $\Eis^{gk}(\beta)$ instead of $\Eis^k(\beta)$.
If we consider $e^*\pol^{D}_{\beta}$ as a 
power series in the ${\bf e}_i$ 
we are interested in the coefficient of 
$\frac{N{\bf e}^{\otimes k}}{k!^g}$. 
In fact, the integrallity
properties of  $\Eis^{gk}(\beta)$ are better reflected if we write
it in terms of a basis $a_1,\ldots, a_g$ of $\afr$. Then
$N{\bf e}^{\otimes k}=N\afr^{-k} N{\bf a}^{\otimes k}$, where
${\bf a}_1,\ldots, {\bf a}_g$ denote again the images of 
$a_1,\ldots, a_g$ in $\Uhat(\afr)$.
We get: 
\begin{cor}\label{Eisintegration}
With the above basis ${\bf a}_1,\ldots, {\bf a}_g$, 
The integral over the Eisenstein class is given by
\[
u_{\id,*}\Eis^{gk}(\beta)=
\frac{(-1)^{g-1}(k!)^g}{(2\pi i)^{g(k+1)}|d_F|^{1/2}N\afr^{k+1}}
\sum_{\rho\in \Oh^*_{(n)}\backslash\bigl(\afr^\lor\setminus 0\bigr)}
\frac{e^{-2\pi i\rho(\beta)}}
{N(\rho)^{k+1}}\frac{N{\bf a}^{\otimes k}}{k!^g}.
\]
\end{cor}
\subsection{5. Step: End of the proof}
To finish the proof of the theorem \ref{maintheorem}, 
let $\alpha\in L[\Ah[n]\setminus e(S)]^0$
and suppose we want to compute $\res(\Eis^k(\alpha))(h)$. Using the 
equivariance of $\res\circ \Eis^k$ from (\ref{resequiv}), this amounts to
compute $\res(\Eis^k(h\alpha))(\id)$. Theorem \ref{degenercompu} shows
that 
\[
\res(\Eis^k(h\alpha))(\id)=u_{\id,*}q_*\Eis^k(h\alpha),
\]
where $q:S^1_B\to S^1_T$ and 
$u_{\id}:\Gamma_T\backslash \bigl(T^1(\R)/K^{T^1}_\infty\bigr) \to pt$ is 
the structure map of the component corresponding to $\id\in T^1(\Z/n\Z)$.
From theorem \ref{eisendeg} we get
\[
q_*\Eis^k(h\alpha)=\Eis^k(p(h\alpha)).
\]
Using corollary \ref{Eisintegration} for $\afr=\Oh$
and the formula \ref{fcteq} for $\bfr=\ffr=\Oh$ we get
\begin{equation}\label{finalformula}
\frac{(-1)^{g-1}(k!)^g}{(2\pi i)^{gk+g}|d_F|^{1/2}}
\sum_{\rho\in \Oh^*_{(n)}\backslash\bigl(\Oh^\lor\setminus 0\bigr)}
\frac{e^{-2\pi i\rho(p(h\alpha))}}
{N(\rho)^{k+1}}=(-1)^{g-1}\sum_{\sigma\in D}l_\sigma \zeta(\Oh,\Oh,p(h\sigma),-k),
\end{equation}
which is the formula in the main theorem \ref{maintheorem}.
To prove the corollary, we use that the map of real tori
\[
\Ah(\C)\xrightarrow{p}\Th_\Mh
\]
factors through $\phi:\Th_{\bfr_\tildeh}\to\Th_\Mh $, where $\phi$ is
induced by the inclusion $\bfr_\tildeh\subset \Oh$.
Using corollary \ref{Eisintegration} for $\afr=\bfr_\tildeh$, we get the
desired formula
\[
\res(\Eis^{gk}(\alpha))(h)=(-1)^{g-1}N\bfr_\tildeh^{-k-1}
\sum_{\sigma\in D}l_\sigma \zeta(\bfr_\tildeh,\Oh,p_\tildeh(\sigma),-k),
\]
which ends the proof.

\end{document}